\documentclass{llncs}

\long\def\ignore#1{\relax}

\usepackage{a4}
\usepackage{amssymb}
\usepackage{QED}
\usepackage{url}

\ignore{

\newtheorem{theorem}{Theorem}

\newtheorem{definition}[theorem]{Definition}

\newtheorem{lemma}[theorem]{Lemma}

\newtheorem{proposition}[theorem]{Proposition}

}

\def\be{\beta}

\def\l{\lambda}

\def\G{\Gamma}

\def\f{\rightarrow}

\def\v{\vdash}
\def\<{\langle}
\def\>{\rangle}
\def\F{\displaystyle\frac}

\def\R{\ifmmode{\rm I\mkern-3.1mu R\mkern1mu}\else{\rm
I\kern-.18em  R\hskip1pt\ }\fi\relax}

\def\Z{\ifmmode{ Z\mkern-8.0mu Z\mkern2mu}\else{
Z\kern-.32em Z\hskip1pt\ }\fi\relax}

\def\Q{\ifmmode{\rm Q\mkern-10mu
l\mkern4.5mu}\else{\rm Q\kern-.57em l\hskip3pt\
}\fi\relax}

\def\N{\ifmmode{\rm I\mkern-3.1mu
N\mkern0.5mu}\else{\rm I\kern-.16em N\hskip0.5pt\
}\fi\relax}

\def\C{\ifmmode{\rm C\mkern-8.8mu l\mkern4mu}\else{\rm
C\kern-.48em l\hskip2.6pt\ }\fi\relax}

\def\m{\mu}
\def\al{\alpha}
\def\tr{\triangleright}
\def\la{\lambda}
\def\ras{\tr^*}
\def\ih{{\em IH}}
\def\ihb{{\em IH }}

\title{An arithmetical proof of the strong
normalization  \\for the  $\l$-calculus \\ with recursive
equations on types}

\author{Ren\'e David \&
Karim Nour\thanks{Universit\'e de Savoie, Campus Scientifique,
73376 Le Bourget du Lac, France. \hspace{2cm}
 Email : \{david,
nour\}@univ-savoie.fr}}

\institute{Universit\'e de Savoie}

\begin{document}
\maketitle \pagestyle{plain}

\begin{abstract}
We give an arithmetical proof of the strong normalization of the
$\l$-calculus (and also of the $\l\m$-calculus) where the type
system is the one of simple types with recursive equations on
types.

The proof  using candidates of reducibility is an easy extension
of the one without equations but this proof cannot be formalized
in Peano arithmetic. The strength of the system needed for such a
proof was not known. Our proof shows that it is not more than
Peano arithmetic.
\end{abstract}

\section{Introduction}

The  $\l$-calculus is a powerful model for representing functions.
In its un-typed version, every recursive function can be
represented. But, in this model,  a term can be applied to itself
and   a computation may not terminate. To avoid this problem,
types are used. In the simplest case, they are built from atomic
types with the arrow and the typing rules say that a function of
type $U \f V$ may only be applied to an argument of  type $U$.
This discipline ensures that every typed term is strongly
normalizing, i.e. a computation always terminate.

In this system (the simply typed $\l$-calculus), Church numerals,
i.e. the terms of the form $\l f \l x (f \; (f \; ... \; (f \;
x)))$, are codes for the integers. They are the only terms (in
normal form) of type  $(o \f o) \f (o \f o)$. Thus, functions on
the integers can be represented but Schwichtenberg \cite{Schw} has
shown  that very few functions are so. He showed that the extended
polynomials (i.e. polynomials with positive coefficients together
with a conditional operator) are the only functions that can be
represented there. Other type systems were then designed  to allow
the representation of more functions. They are built in different
ways.

The first one consists in extending the set of terms. For example,
in G\"{o}del system $T$,  the terms use the usual constructions of
the $\l$-calculus, the constant $0$, the constructor $S$  and an
operator  for recursion. The types are built from the atomic type
$N$ with the arrow.
  This system
represents exactly the functions whose totality can be shown in
Peano first order arithmetic.

The second one consists in keeping the same terms but extending
the type system. This is, for example,  the case of Girard system
$F$ where the types can use a second order universal quantifier.
There, the type of the integers is given by $\forall X \, ((X \f
X) \f (X \f X))$. This system represents exactly the functions
whose totality can be shown in Peano second order arithmetic.

A third way consists in extending the {\em logic}. In the
Curry-Howard correspondence, the previous systems correspond to
{\em intuitionistic} logic. Other systems correspond to {\em
classical} logic. There, again, new constructors for terms are
introduced. This is, for example, the case of Parigot's
$\l\m$-calculus \cite{Par4}.

Since the introduction of Girard system $F$ for intuitionistic
logic and Parigot's $\l\m$-calculus for classical logic, many
others, more and more powerful, type systems were introduced. For
example, the calculus of constructions (Coquand \& Huet
\cite{Coq}) and, more generally,  the Pure Type Systems.

It is also  worth here to mention the system $TTR$ of Parigot
\cite{Par3} where some types are defined as the least fixed point
of an operator. This system was introduced, not to represent more
functions, but to represent more {\em algorithms}. For example, to
be able to represent the integers in such a way that the
predecessor can be computed in constant time, which is not the
case for the previous systems.

These systems all satisfy the subject reduction (i.e. the fact
that the type is preserved by reduction), the strong normalization
(i.e. every computation terminates) and, for the systems based on
simple types, the decidability of type assignment.

We study here  other kinds of extension of the simply typed
$\l$-calculus, i.e. systems where {\em equations} on types are
allowed. These types are  usually called {\em recursive types}.
For more details see, for example, \cite{Bar3}. They are present
in many languages and are intended to be able to be {\em unfolded}
recursively to match other types. The subject reduction and the
decidability of type assignment are preserved but the strong
normalization may be lost. For example, with the equation $X = X
\f T$, the term $(\delta \; \delta)$ where $\delta=\l x \;(x \,
x)$  is typable but is not strongly normalizing. With the equation
$X=X\f X$, every term can be typed.

By making some natural assumptions on the recursive equations the
strong normalization can be preserved. The simplest condition is to
accept the equation $X = F$ (where $F$ is a type containing the
variable $X$) only when the variable $X$ is positive in $F$. For a set
$\{X_i = F_i \;/ \; i \in I\}$ of mutually recursive equations,
Mendler \cite{Mend1} has given a very simple and natural condition
that ensures the strong normalization of the system.  He also showed
that the given condition is necessary to have the strong
normalization.  His proof is based on the reducibility method. The
condition ensures enough monotonicity to have fixed point on the
candidates. But this proof (using candidates of reducibility) cannot
be formalized in Peano arithmetic and the strength of the system
needed for a proof of the strong normalization of such systems was not
known.

In this paper, we give an {\em arithmetical} proof of the strong
normalization of the simply typed $\l$-calculus (and also of the
$\l\mu$-calculus) with recursive equations on types satisfying
Mendler's condition.

This proof is  an extension of the one given by the first author
for the simply typed $\l$-calculus. It can be found either in
\cite{Dav1} (where it appears among many other things) or as a
simple unpublished note on the web page of the first author
\cite{Dav2}. Apparently, proof methods similar to that used here
were independently invented by several authors (Levy, van Daalen,
Valentini and others). The proof for the $\l\m$-calculus is an
extension of the ones given in \cite{DN2} or \cite{DN3}.

The paper is organized as follows. In section \ref{S-lambda} we
define the simply typed $\l$-calculus with recursive equations on
types. To help the reader and show the main ideas, we first give,
in section \ref{proof}, the proof of strong normalization for the
$\l$-calculus. We generalize this proof to the $\l\mu$-calculus in
section \ref{S-mu}. In section \ref{applic}, we give two examples
of applications of systems with recursive types.  We conclude in
section \ref{conc} with some open questions.

\section{The typed $\l$-calculus}\label{S-lambda}

\begin{definition}
Let ${\cal V}$ be an infinite set of variables.
\begin{enumerate}
\item The set ${\cal M}$ of $\l$-terms is defined
by the following grammar

$${\cal M} ::= \; {\cal V} \; \mid \; \l {\cal V} \
{\cal M} \; \mid
\; ({\cal M} \; {\cal M})$$
\item The relation $\rhd$ on ${\cal M}$ is defined
as the least  relation (compatible with the context) containing
the rule $(\l x \; M \; N) \rhd M[x:=N]$. As usual, $\rhd^*$
(resp. $\tr^+$)  denotes the reflexive and transitive (resp.
transitive) closure of $\rhd$.
\end{enumerate}
\end{definition}

\begin{definition}
Let ${\cal A}$ be a set of atomic constants and ${\cal X}=\{X_i \;
/ \; i \in I \}$ be a set of type variables.
\begin{enumerate}
\item The set ${\cal T}$ of types is defined by the following grammar

$${\cal T} ::= \; {\cal A} \; \mid \; {\cal X} \; \mid {\cal T} \f
{\cal T}$$

\item When $E=\{F_i \; /
\; i \in I \}$ is a set of types,  the congruence $\approx$
generated by $E$ is the least congruence on ${\cal T}$ such that
 $X_i \approx F_i$ for each $i \in I$.

\end{enumerate}
\end{definition}

\begin{definition}

Let $\approx$ be a congruence on ${\cal T}$. The typing rules of
the typed system  are given below where $\G$ is a context, i.e. a
set of declarations of the form $x : U$ where $x \in {\cal V}$ and
$U \in {\cal T}$.

\begin{center}
$\F{}{\G , x : U \v x : U}\;\;\;ax$ \hspace{1cm}  $\F{\G \v M: U
\;\;\; \hspace{0.2in} U \approx V}{\G \v M : V}\;\;\;\approx$
\medskip

$\F{\G , x : U \v M : V}{\G \v \l x \; M : U \f
V}\;\;\;\f_i$\hspace{1cm} $\F{\G \v M_1 : U \f V \;\;\;
\hspace{0.2in} \G \v M_2 : U}{\G \v (M_1 \; M_2) : V}\;\;\;\f_e$
\end{center}

\end{definition}
\begin{lemma}\label{congru}
Let $\approx$ be a congruence generated by a set of types.
\begin{enumerate}
\item If $U \approx V_1 \f V_2$, then $U \in {\cal X}$ or $U = U_1 \f
U_2$.
\item If $U_1 \f V_1 \approx U_2 \f V_ 2$, then $U_1 \approx U_2$ and
$V_1 \approx V_2$.
\item If $\G \v x: T$, then $x : U$ occurs in $\G$ for
some $U \approx T$.
\item If $\G \v \l x \; M : T$, then  $\G , x : U \v M :
V$ for some $U,V$ such that $U \f V \approx T$.
\item If $\G \v (M \, N) : T$, then $\G  \v M : U \f
V$, $\G  \v N : U$ for some $V  \approx T$ and $U$.
\item If $\G , x : U\v M: T$ and $U \approx V$, then
$\G , x : V \v M: T$.
\item If $\G , x : U \v M : T$ and $\G \v N : U$, then
 $\G  \v
M[x:= N] : T$.
\end{enumerate}
\end{lemma}

\begin{proof}
Easy.
\end{proof}

\begin{theorem}\label{SR}
If $\G \v M:T$ and $M \rhd^* M'$, then $\G \v M':T$.
\end{theorem}

\begin{proof}
It is enough to show that if $\G \v (\l x \; M \, N):T$, then $\G
\v M[ x:= N]:T$. Assume $\G \v (\l x \; M \, N):T$. By lemma
\ref{congru}, $\G \v \l x \; M : U \f V$, $\G \v N:U$ and $V
\approx T$. Thus, $\G , x : U' \v M : V'$ and $U' \f V' \approx U
\f V$. By lemma \ref{congru}, we have $U' \approx U$ and $V'
\approx V$. Thus,  $\G , x : U \v M : V$. Since $\G \v N:U$ and $V
\approx T$, the result follows immediately.
\end{proof}

\begin{definition}
Let $X \in {\cal X}$. We define the subsets ${\cal T}^+(X)$ and
${\cal T}^-(X)$ of ${\cal T}$ as follows.
\begin{itemize}
\item $X \in {\cal T}^+(X)$
\item
If $\; U \in ({\cal X}- \{X\}) \cup \cal{A}$, then $U \in {\cal
T}^+(X) \cap {\cal T}^-(X)$.
\item If $U  \in {\cal T}^-(X)$ and $V \in {\cal
T}^+(X)$, then $U \f V \in {\cal T}^+(X)$ and  $V \f U \in {\cal
T}^-(X)$.
\end{itemize}
\end{definition}

\begin{definition}\label{cong}
 We say that a congruence $\approx$ is
good if the following property holds: for each $X \in {\cal X}$,
if $X \approx T$, then $T \in {\cal T}^+(X)$.
\end{definition}

\noindent {\bf Examples}

In each of the following cases, the congruence generated by the given equations is good.

\begin{enumerate}

\item $X_1 \approx (X_1 \f X_2 \f Y) \f Y$ and $X_2 \approx (X_2 \f X_1 \f Y) \f Y$.

\item $X_1 \approx X_2 \f X_1$ and $X_2 \approx X_1 \f X_2$.

\item The same equations as in case 2 and $X_3 \approx F(X_1,X_2) \f
  X_3$ where $F$ is any type using only the variables $X_1,X_2$.

\item The same equations as in case 3 and $X_4 \approx X_5 \f
G(X_1,X_2,X_3) \f X_4$, $X_5 \approx X_4 \f H(X_1,X_2,X_3) \f X_5$
where $G,H$ are any types using only the variables $X_1,X_2,X_3$.

\end{enumerate}

{\em In the rest of the paper, we fix a finite set $E=\{F_i \; /
\; i \in I \}$ of types and we denote by $\approx$ the congruence
generated by $E$.  We assume that $\approx$ is good.}

\bigskip

\noindent {\bf Notations and remarks}
\begin{itemize}

 \item We  have assumed that the set of equations that we consider
is finite. This is to ensure that the order on $I$ given by
definition \ref{ordre} below is well founded. It should be clear
that this is not a real constraint. Since to type a term, only a
finite number of equations is used, we may  consider that the
other variables are constant and thus the general result follows
immediately from the finite case.

\item If $M$ is a term, $cxty(M)$ will denote the structural
complexity of $M$.

\item We denote by $SN$ the set of  strongly
normalizing terms. If $M \in SN$, we denote by  $\eta(M)$ the
length of the longest reduction of $M$ and by $\eta c(M)$ the pair
$\langle \eta(M), cxty(M)\rangle$.
\item We denote by $M \preceq N$ the fact that $M$ is a sub-term of a
reduct of $N$.
\item As usual, some parentheses are omitted and, for example,  we write $(M \; P \; Q)$ instead
of $((M \; P) \; Q)$. More generally, if $\overrightarrow{O}$ is a
finite sequence $O_1,...,O_n$ of terms, we denote by $(M \;
\overrightarrow{O})$ the term $((...(M \; O_1)  ... \; O_{n-1} )\;
O_n)$ and by $\overrightarrow{O} \in SN$ the fact that $O_1, ...,
O_n \in SN$.

\item If $\sigma$ is the substitution $[x_1:=N_1, ...,
x_n:=N_n]$, we denote by $dom(\sigma)$ the set $\{x_1,...,x_n\}$,
by $Im(\sigma)$ the set $\{N_1, ..., N_n\}$ and by $\sigma \in SN$
the fact that $Im(\sigma) \subset SN$.

\item If $\sigma$ is a substitution, $z \not \in dom(\sigma)$
and $M$ is a term, we denote by $[\sigma +
z:=M]$ the substitution $\sigma'$ defined by
$\sigma'(x)=\sigma(x)$ for $x \in dom(\sigma)$ and
$\sigma'(z)=M$.
\item In a proof by induction, \ihb will denote the
induction hypothesis. When the induction is done on a tuple of
integers, the order always is the lexicographic order.

\end{itemize}

\section{Proof of the strong normalization}\label{proof}

\subsection{The idea of the proof}\label{3.1}

We give the idea for one equation $X \approx F$.  The extension
for the general case is given at the beginning of section
\ref{3.4}.

It is enough to show that, if $M,N$ are in $SN$, then $M[x:=N] \in
SN$. Assuming it is not the case, the interesting case is $M= (x
\; P )$ with $(N \; P_1) \not \in SN$ where $P_1 =P[x:=N] \in SN$.
This implies that $N \ras \l y N_1$ and $N_1[y=P_1] \not \in SN$.
If we know that the type of $N$ is an arrow type, we get a similar
situation to the one we started with, but where the type of the
substituted variable has decreased.  Repeating the same argument,
we get the desired result, at least for $N$ whose type does not
contain $X$. If it is not the case, since, by repeating the same
argument, we cannot come to a constant type (because such a term
cannot be applied to something), we come to $X$. Thus, it remains
to show that, if $M,N$ are in $SN$ and the type of $x$ is $X$,
then $M[x:=N] \in SN$.

To prove this, we prove something a bit more general. We prove
that, if $M,\sigma \in SN$ where $\sigma$ is a substitution such
that the types of its image are in ${\cal T}^+(X)$, then
$M[\sigma] \in SN$. The proof is done, by induction on $ \eta
c(M)$ as follows. As before, the interesting case is $M=(x \; P),
\sigma(x)=N \ras \l y N_1$, $P_1=P[\sigma] \in SN$ and $N_1[y=P_1]
\not\in SN$. Thus, there is a sub-term of a reduct of $N_1$ of the
form $(y \; N_2)$ such that $(P_1 \; N_2[y:=P_1 ]) \not \in SN$
but $N_2[y:=P_1 ] \in SN$.  Thus $P_1$ must reduce to a~$\l$.

This $\l$ cannot come from some $x' \in dom(\sigma)$, i.e. $ P \;
\ras (x' \; \overrightarrow{Q})$.  Otherwise, the type of $P$
would be both positive (since $ P \ras (x' \; \overrightarrow{Q})$
and the type of $x'$ is positive) and negative (since, in $M$, $P$
is an argument of $x$ whose type also is positive). Thus the type
of $P_1$ (the same as the one of $P$) does not contain  $X$. But
since $N_1, P_1$ are in $SN$, we already know that $N_1[y=P_1]$
must be in $SN$. A contradiction.  Thus, $P \ras \l x_1 M_1$ and
we get a contradiction from the induction hypothesis since we have
$M_1[\sigma'] \not \in SN$ for $M_1$ strictly less than $M$. The
case when $y$ has more than one argument is intuitively treated by
``repeat the same argument'' or, more formally, by lemma
\ref{prepa} below.

As a final remark, note that many lemmas are stated in a negative
style and thus may seem to hold only classically. This has been
done in this way because we believe that this presentation is
closer to the intuition. However, it is not difficult to check
that the whole proof can be presented and done in a constructive
way.

\subsection{Some useful lemmas on the un-typed calculus}\label{3.2}

\begin{lemma}\label{untyped}
Assume $M,N,\overrightarrow{O} \in SN$ and
$(M\;N\;\overrightarrow{O}) \not \in SN$. Then, for some term
$M'$, $M \rhd^* \lambda x \; M'$ and
$(M'[x:=N]\;\overrightarrow{O}) \not \in SN$.
\end{lemma}

\begin{proof}
Since $M,N,\overrightarrow{O} \in SN$, an infinite reduction of $P
= (M\;N\;\overrightarrow{O})$ looks like $P \rhd^* (\lambda x \;
M'\; N'  \; \overrightarrow{O'}) \rhd (M'[x:=N']\;
\overrightarrow{O'}) \rhd \ldots $ and the result immediately
follows from the fact that $(M'[x:=N]\;\overrightarrow{O} ) \rhd^*
(M'[x:=N']\;\overrightarrow{O'})$.
\end{proof}

\begin{lemma}\label{subterm}
Let $M$ be a term and $\sigma$ be a substitution.
Assume $M,
\sigma \in SN$ and $M[\sigma] \not\in SN$. Then
$(\sigma(x)
\;\overrightarrow{P[\sigma]}) \not \in SN$ for some
$(x \;
\overrightarrow{P}) \preceq M$ such that
$\overrightarrow{P[\sigma]} \in SN$.
\end{lemma}

\begin{proof}
A sub-term $M'$ of a reduct of $M$ such that $\eta c(M')$ is
minimum  and $M'[\sigma] \not\in SN$ has the desired form.
\end{proof}

\begin{lemma}\label{lambda}
Let $M$ be a term and $\sigma$ be a substitution such that
$M[\sigma] \rhd^* \lambda z M_1$. Then \\
- either $M \rhd^* \lambda
z M_2$ and $M_2[\sigma] \rhd^* M_1$\\
- or  $M \rhd^* (x \; \overrightarrow{N})$ for some $x \in
dom(\sigma)$ and $(\sigma(x) \; \overrightarrow{N[\sigma]}) \rhd^*
\lambda z M_1$.
\end{lemma}

\begin{proof}
This is a classical (though not completely trivial) result in
$\l$-calculus. Note that, in case $M \in SN$ (and we will only use
the lemma in this case), it becomes easier. The proof can be done
by induction on $\eta c(M)$ by considering the possibility for
$M$: either $\l y M_1$ or $(\l y M_1 \; P \; \overrightarrow{Q})$
or $(x \; \overrightarrow{N})$ (for $x$ in
 $ dom(\sigma)$ or not).
\end{proof}

\subsection{Some useful lemmas on the congruence}\label{3.3}

\begin{definition}\label{ordre}
We define on $I$ the following relations
\begin{itemize}
\item $i \leq j$ iff $X_i \in var(T)$ for some $T$ such that $X_j \approx
  T$.
\item $i \sim j$ iff $i \leq j$ and $j \leq i$.
 \item $i < j$ iff $i \leq j$ and $j \not\sim i$
\end{itemize}
\end{definition}

It is clear that $\sim$ is an equivalence on $I$.


\begin{definition}
\begin{enumerate}
\item Let ${\cal X}_i= \{X_j \;/ \; j \leq i\}$ and
${\cal X}'_i= \{X_j \;/ \; j < i\}$.
\item For ${\cal Y} \subseteq {\cal X}$, let ${\cal T}({\cal Y}) =\{T \in {\cal T} \; / \; var(T) \subseteq {\cal Y}\}$
 where $var(T)$ is the set of type
variables occurring in $T$.
 \item For $i\in I$, we will abbreviate by ${\cal T}_i$ the set ${\cal
  T}({\cal X}_i)$ and by ${\cal T}'_i$ the set ${\cal T}({\cal X}'_i)$.
  \item If $\varepsilon \in \{+,-\}$,
$\overline{\varepsilon}$ will  denote the opposite of
$\varepsilon$. The opposite of + is - and conversely.
\end{enumerate}

\end{definition}

\begin{lemma}
Let $i \in I$.  The class of $i$ can be partitioned into two
disjoint sets $i^+$ and $i^-$ satisfying the following properties.
\begin{enumerate}
\item If $\varepsilon \in \{+,-\}$, $j \in i^{\varepsilon}$ and $X_j
\approx T$, then for each $k \in i^{\varepsilon}$, $T \in {\cal
T}^{\varepsilon}(X_k)$ and for each $k \in
i^{\overline{\varepsilon}}$, $T \in {\cal
T}^{\overline{\varepsilon}}(X_k)$.
\item Let $j\sim i$. Then, if $j \in i^+$, $j^+=i^+$ and $j^-=i^-$ and
if $j \in i^-$, $j^+=i^-$ and $j^-=i^+$.
\end{enumerate}
\end{lemma}

\begin{proof}
This follows immediately from the following observation. Let $i \sim
 j$ and $X_i \approx T \approx U$.  Choose an occurrence of $X_j$ in
 $T$ and in $U$.  Then, these occurrences have the same polarity. This
 is because, otherwise, since $i \leq j$, there is a $V$ such that
 $X_j \approx V$ and $X_i$ occurs in $V$. But then, replacing the
 mentioned occurrences of $X_j$ by $V$ in $T$ and $U$ will contradict
 the fact that $\approx$ is good.
\end{proof}

\begin{definition}\label{class}
Let $i \in I$  and $\varepsilon \in \{+,-\}$. We
denote ${\cal T}_i^{\varepsilon} = \{T \in {\cal T}_i$ / for each $j \in
i^{\varepsilon}$, $T \in {\cal T}^{\varepsilon}(X_j)$ and for each $j
\in i^{\overline{\varepsilon}}$, $T \in {\cal
T}^{\overline{\varepsilon}}(X_j)\}$.
\end{definition}

\begin{lemma}\label{propr}
Let $i \in I$  and $\varepsilon \in \{+,-\}$.
\begin{enumerate}
\item $ {\cal T}_i^{\varepsilon} \cap {\cal
  T}_i^{\overline{\varepsilon}}\subseteq {\cal T}'_i$.
\item If $U \in {\cal T}_i^{\varepsilon}$ and $U \approx
 V$, then $V \in {\cal T}_i^{\varepsilon}$.
\item If $U \in {\cal T}_i^{\varepsilon}$ and $U \approx U_1 \f U_2$,
  then $U_1 \in {\cal T}_i^{\overline{\varepsilon}}$ and $U_2 \in
  {\cal T}_i^{\varepsilon}$.
\end{enumerate}
\end{lemma}

\begin{proof}
Immediate.
\end{proof}

\medskip

\noindent {\bf Notations, remarks and examples}
\begin{itemize}
\item If the equations are those of the case 4 of the examples
given above, we have $1 \sim 2 < 3 < 4 \sim 5$ and, for example,
$1^+=\{1\}$ and $1^-=\{2\}$, $3^+=\{3\}$, $3^-=\emptyset$,
$4^+=\{4\}$ and $4^-=\{5\}$.
\item If $T$ is a type, we denote by $lg(T)$ the size of $T$. Note
that the size of a type is, of course, not preserved by the
congruence. The size of a type will only be used in lemma
\ref{sansX} and the only property  that we will use  is that
$lg(U_1)$ and $lg(U_2)$ are less than $lg(U_1 \f U_2)$.
\item By the typing rules, the type of a term can be freely replaced by
an equivalent one. However, for $i \in I$ and $\varepsilon \in
\{+,-\}$, the fact that $U \in {\cal T}_i^{\varepsilon}$ does not
change when $U$ is replaced by $V$ for some $V \approx U$. This
will be used extensively in the proofs of the next sections.

\end{itemize}

\subsection{Proof of the strong normalization}\label{3.4}
To give the idea of the proof, we first need a definition.

\begin{definition}
Let $\cal{E}$ be  a set of types. Denote by $H[\cal{E}]$ the
following property:

Let $M,N \in SN$. Assume $\G, x:U \v M : V$ and $\G \v N : U$ for
some $\G, U, V$ such that $U \in \cal{E}$. Then $M[x:=N] \in SN$.
\end{definition}

To get the result, it is enough to show $H[{\cal T}]$. The proof that
 any typed term is in $SN$ is then done by induction on $cxty(M)$. The
 only non trivial case is $M=(M_1 \; M_2)$. But $M= (x \;
 M_2)[x:=M_1]$ and the result follows from $H[{\cal T}]$ and the {\em
 IH}.

\medskip

We first show the following (see lemma \ref{sansX}). Let ${\cal Y}
\subseteq {\cal X}$. To prove $H[{\cal T}({\cal Y})]$, it is enough to
prove $H[\{X\}]$ for each $X \in {\cal Y}$.

\medskip

It is thus enough to prove of $H[\{X_i\}]$ for each $i \in I$.
This is done by induction on $i$. Assume $H[\{X_j\}]$ for each $j
< i$.  Thus, by the previous property, we know $H[{\cal T}'_i]$.
We show $H[\{X_i\}]$ essentially as we said in section \ref{3.1}.
The only difference is that, what was called there `` $X$ is both
positive and negative in $T$'' here means $T$ is both in ${\cal
T}_i^+$ and ${\cal T}_i^-$. There we deduced that $X$ does not
occur in $T$. Here we deduce $T \in {\cal T}'_i$ and we are done
since we know the result for this set.

\begin{lemma}\label{sansX}
 Let ${\cal Y} \subseteq {\cal X}$ be such that $H[\{X\}]$ holds for
 each $X \in {\cal Y}$. Then $H[{\cal T}({\cal Y})]$ holds.
\end{lemma}

\begin{proof}
Let $M,N$ be terms in $SN$. Assume $\G, x:U \v M : V$ and $\G \v N :
U$ and $U \in {\cal T}({\cal Y})$.  We have to show $M[x := N] \in SN$.

This is done by induction on $lg(U)$. Assume $M[x:=N] \not\in SN$.  By
  lemma \ref{subterm}, let $(x \, P \, \overrightarrow{Q})\preceq M$
  be such that $P_1, \overrightarrow{Q_1} \in SN$ and $(N \, P_1 \,
  \overrightarrow{Q_1}) \not \in SN$ where $P_1=P[x:=N]$ and
  $\overrightarrow{Q_1} = \overrightarrow{Q[x:=N]}$.  By lemma
  \ref{untyped}, $N \rhd^* \l x_1 N_1$ and $(N_1[x_1:=P_1] \;\;
  \overrightarrow{Q_1}) \not \in SN$.

If $U$ is a variable (which is in ${\cal Y}$ since $U \in {\cal
T}({\cal Y})$), we get a contradiction since we have assumed that
$H[\{X\}]$ holds for each $X \in {\cal Y}$.

The type $U$ cannot be a constant since, otherwise $x$ could not be
applied to some arguments.

Thus $U = U_1 \f U_2$. In the typing of $(N \, P_1 \,
\overrightarrow{Q_1})$, the congruence may have been used and
thus, by lemma \ref{congru}, there are $W_1 \approx U_1$, $W_2
\approx U_2$, $U \approx W_1 \f W_2$ and $\G, x_1: W_1 \v N_1 :
W_2$ and $\G \v P_1 : W_1$. But then, we also have $\G, x_1: U_1
\v N_1 : U_2$ and $\G \v P_1 : U_1$. Now, by the {\em IH}, we have
$N_1[x_1:=P_1] \in SN$ since $lg(U_1) < lg(U)$. Since $\G, z : U_2
\v (z \; \overrightarrow{Q_1}) : V'$ for some $V'$ and $\G \v
N_1[x_1:=P_1] : U_2$, by the {\em IH} since $lg(U_2) < lg(U)$, we
have $(N_1[x_1:=P_1] \; \overrightarrow{Q_1})=(z \;
\overrightarrow{Q_1})[z=N_1[x_1:=P_1]] \in SN$. Contradiction.
\end{proof}

\medskip

\noindent {\em For now on, we fix some $i$ and we assume
$H[\{X_j\}]$ for each $j < i$.  Thus, by lemma \ref{sansX}, we
know that $H[{\cal T}'_i]$ holds. It remains to prove $H[\{X_i\}]$
i.e. proposition \ref{avecX}.}

\medskip

\begin{definition}
Let $M$ be a term, $\sigma$ be a substitution, $\G$ be a context and
$U$ be a type. Say that $(\sigma,\G , M , U)$ is adequate if the
following holds.
\begin{itemize}
\item  $\G \v M[\sigma] :
U$ and $M, \sigma \in SN$.
\item For each $x \in dom(\sigma)$,   $\G \v
\sigma(x) : V_x$ and $V_x \in {\cal T}_i^+ $.
\end{itemize}
\end{definition}

\begin{lemma}\label{prepa}
Let $n,m$ be integers, $\overrightarrow{S}$ be a sequence of terms
and $(\delta, \Delta , P , B)$ be adequate.  Assume that
\begin{enumerate}
\item  $B \in {\cal T}_i^- -  {\cal T}'_i$ and
$\Delta \v (P [\delta] \, \overrightarrow{S}) : W$ for some $W$.
\item $\overrightarrow{S} \in SN$, $P \in SN$ and $\eta c(P) < \langle n, m \rangle$.
\item $M[\sigma] \in SN$ for every adequate $(\sigma, \G, M, U)$ such
that $\eta c(M) < \langle n, m \rangle  $.
\end{enumerate}
Then $(P [\delta] \; \overrightarrow{S})  \in SN$.
\end{lemma}

\begin{proof}
By induction on the length of $\overrightarrow{S}$. If
$\overrightarrow{S}$ is empty, the result follows from (3) since
$\eta c(P) < \langle n, m \rangle$.  Otherwise, let
$\overrightarrow{S} = S_1 \overrightarrow{S_2}$ and assume that $P
[\delta] \rhd^* \l z \; R$. By lemma \ref{lambda}, there are two
cases to consider:
\begin{itemize}
\item $P \rhd^* \l z \; R'$. We have to show that $Q = (R' [\delta +
z:= S_1]\ \overrightarrow{S_2}) \in SN$. Since $B \in {\cal
T}_i^-$, by lemmas \ref{congru} and \ref{propr}, there are types
$B_1, B_2$ such that $B \approx B_1 \rightarrow B_2$ and $\Delta,
z : B_1 \v R' : B_2 $ and $\Delta \v S_1 : B_1$ and $B_1 \in {\cal
T}_i^+$ and $B_2 \in {\cal T}_i^-$.  Since $\eta c(R') < \langle
n, m \rangle$ and $([\delta +z=S_1], \Delta \cup \{ z : B_1\}, R',
B_2)$ is adequate, it follows from (3) that $R' [\delta + z:=
S_1]] \in SN$.

- Assume first $B_2 \in {\cal T}'_i$. Since $(z' \,
\overrightarrow{S_2}) \in SN$   and $Q = (z' \,
\overrightarrow{S_2}) [z':= R' [\delta + z:= S_1]]$, the result
follows from $H[{\cal T}'_i]$.

-  Otherwise, the result follows from the {\em IH} since $([\delta
+z=S_1], \Delta \cup \{ z : B_1\}, R', B_2)$ is adequate and the
length of $ \overrightarrow{S_2}$ is less than the one of
$\overrightarrow{S}$.

\item If $P \rhd^* (y \, \overrightarrow{T})$ for some $y \in
dom(\delta)$. Then $\Delta \v (\delta(y) \,
\overrightarrow{T[\delta]}) : B $. By the definition of adequacy,
the type of $y$ is in $ {\cal T}_i^+$ and   $B \in {\cal T}_i^-
\cap {\cal T}_i^+ \subseteq {\cal T}'_i$. Contradiction. \qed
\end{itemize}
\end{proof}

\begin{lemma}\label{X>0}
Assume $(\sigma,  \G, M,A )$ is adequate.  Then $M
[\sigma] \in
SN$.
\end{lemma}

\begin{proof}
By induction on $\eta c(M)$. The only non trivial case is $M =
(x\; Q\;\overrightarrow{O})$ for some $x \in dom(\sigma) $. Let
$N=\sigma(x)$.

By the {\em IH}, $Q [\sigma] , \overrightarrow{O[\sigma]} \in SN$.
By lemma \ref{congru}, we have $V_x \approx W_1 \f W_2$, $\G \v
Q[\sigma] : W_1$ and $\G \v (N \; Q[\sigma]) : W_2$. Moreover, by
lemma \ref{propr}, $W_1 \in {\cal T}_i^-$ and $W_2 \in {\cal
T}_i^+$. Since $M [\sigma] = (z \; \overrightarrow{O})[\sigma + z
: = (N \, Q[\sigma])]$, $\eta((z \; \overrightarrow{O})) \leq
\eta(M)$, $cxty((z \; \overrightarrow{O})) < cxty(M)$ and $W_2 \in
{\cal T}_i^+$, it is enough, by the {\em IH}, to show that $(N \,
Q[\sigma]) \in SN$.  Assume that $N \rhd^* \lambda y \; N'$. We
have to show that $N'[y:=Q [\sigma]]\in SN$.

- Assume first $W_1 \in {\cal T}'_i$. The result follows from $H[{\cal
  T}'_i]$.

- Otherwise, assume $N'[y:=Q [\sigma]] \not \in SN$. Since $N', Q
[\sigma] \in SN$, by lemma \ref{subterm}, $(y \,
\overrightarrow{L}) \preceq N'$ for some $\overrightarrow{L}$ such
that $ \overrightarrow{L[y:=Q [\sigma]]} \in SN$ and  $(Q [\sigma]
 \; \overrightarrow{L[y:=Q [\sigma]]}) \not \in
SN$. But this contradicts lemma \ref{prepa}. Note that,  by the
\ih, condition (3) of this lemma is satisfied.
\end{proof}

\begin{proposition}\label{avecX}
Assume $\G , x : X_i \v M : U$ and $\G \v N : X_i$ and $M,N \in
SN$. Then $M[x := N] \in SN$.
\end{proposition}

\begin{proof}
This follows from lemma \ref{X>0} since $([x:=N], \G, M, U )$ is
adequate.
\end{proof}

\section{The typed $\lambda\mu$-calculus}\label{S-mu}

\begin{definition}
\begin{enumerate}

\item Let ${\cal W}$ be an infinite set of variables such that
${\cal V} \cap {\cal W} = \emptyset$. An element of ${\cal V}$
(resp. ${\cal W}$) is said to be a $\l$-variable (resp. a
$\mu$-variable). We extend the set of terms by the following rules

$$
{\cal M} ::= ... \mid \mu {\cal W} \, {\cal M} \mid   ({\cal W} \;
{\cal M})
$$

\item We add to the set ${\cal A}$ the constant symbol $\bot$ and we
  denote by $\neg U$ the type $U \f \bot$.

\item  We extend the typing rules by

\begin{center}

$\F{\G , \al : \neg U  \v M : \bot} {\G \v \mu \al M : U }  \,
\bot_e$ \hspace{0.5cm} $\F{\G , \al : \neg U  \v M : U} {\G, \al :
\neg U  \v ( \al \; M) : \bot } \, \bot_i$

\end{center}
\noindent where $\G$ is now a set of declarations of the form $x :
U$ and $\al : \neg U$ where $x$ is a $\l$-variable and $\al$ is a
$\m$-variable.

\item We add to $\tr$ the following reduction rule
$(\mu \al M \; N) \tr \mu \al M[\al= N]$ where $M[\al= N]$ is
obtained by replacing each sub-term of $M$ of the form $(\al \;
P)$ by $(\al \; (P \; N))$. This substitution will be called a
$\m$-substitution whereas the (usual) substitution $M[x:=N]$ will
be called a $\l$-substitution.
\end{enumerate}
\end{definition}

\noindent {\bf Remarks}
\begin{itemize}
\item Note that we adopt here a more liberal syntax
(also called de
Groote's calculus \cite{DeGr}) than in the original calculus since
we do not
ask that a $\m \al$ is immediately followed by a $(\be
\; M)$
(denoted $[\be] M$ in Parigot's notation).

\item We also have changed Parigot's typing notations. Instead of
writing $M : (A_1^{x_1 }, ...$, $ A_n^{x_n } \vdash B,
C_1^{\al_1}, ..., C_m^{\al_m })$ we have written $x_1 : A_1, ...,
x_n : A_n, \al_1: \neg C_1, ..., \al_m : \neg C_m \vdash M : B$
but, since the first introduction of the $\l\m$-calculus, this is
now quite common.

\item Unlike for a $\l$-substitution where, in
$M[x:=N]$, the
variable $x$ has disappeared it
 is important to note that, in a $\m$-substitution,
the variable $\al$ has not disappeared. Moreover its
type has
 changed.  If the type of $N$ is $U$ and, in $ M $,
the type of $\al$
 is $\neg (U\rightarrow V)$ it becomes $\neg V$ in
$M[\al=N]$.
\item The definition of good congruence is the same as before. As a
consequence, we now have the following facts. If $U \approx \bot$,
then $U=\bot$ and, if $\neg U \approx \neg V$, then $ U \approx
 V$.
 \item We also extend all the notations given in section \ref{S-lambda}.
Finally note that lemma \ref{congru}  remains valid. Moreover,
they are easily extended by lemma \ref{typed1'} below.
\end{itemize}

\begin{lemma}\label{typed1'}
\begin{enumerate}
\item If $\G \v \mu \al \; M : U$, then $\G , \al : \neg V \v M : \bot$ for some
$V$ such that $U \approx V$.
\item If $\G , \al : \neg U \v (\al \, M) : T$, then $\G , \al : \neg U
  \v M : U$ and $T = \bot$.
  \item If $\G , \al : \neg (U \f V) \v M : T$ and $\G \v N : U$, then
 $\G , \al : \neg V \v M[\al= N] : T$.
\end{enumerate}
\end{lemma}

\begin{theorem}\label{SR'}
If $\G \v M:T$ and $M \rhd^* M'$, then $\G \v M':T$.
\end{theorem}

\begin{proof}
It is enough to show that, if $\G \v (\mu \al \; M \, N):T$, then
$\G \v \mu \al \, M[ \al= N]:T$. Assume $\G \v (\mu \al \; M \,
N):T$. By lemma \ref{congru}, $\G \v \mu \al \; M : U \f V$, $\G
\v N:U$ and $V \approx T$. Thus,  $\G , \al : \neg T' \v M : \bot$
and $T' \approx U \f V$. By lemma \ref{congru}, we have $\G , \al
: \neg (U \f V) \v M : \bot$. Since $\G \v N:U$ and $V \approx T$,
 $\G , \al : \neg V \v M[\al = N] : \bot$.
Then $\G \v \mu \al \, M[\al = N] : V$ and $\G  \v \mu \al \,
M[\al = N] : T$.
\end{proof}

\subsection{Some useful lemmas on the un-typed calculus}\label{4.1}

\ignore{
\begin{lemma}\label{ll2}
If $(M \; N) \ras \m \al P$, then  $M \ras \m \al M_1$
and
$M_1[\al = N] \ras P$.
\end{lemma}

\begin{proof}
By induction  on the length of the reduction $(M \; N)
\ras \m \al
P$.
\end{proof}
}

\begin{lemma}\label{mu.lambda}
Let $M$ be a term and $\sigma =\sigma_1 \cup \sigma_2$ where
$\sigma_1$ (resp. $\sigma_2$) is $\l$ (resp. $\m$) substitution.
 Assume $M[\sigma] \rhd^* \m\al M_1$ (resp. $\l y M_1$). Then

 - either $M \rhd^* \m
\al M_2$ (resp. $\l y M_2$) and $M_2[\sigma] \rhd^* M_1$

- or ($M \rhd^* (x \; \overrightarrow{N})$ for some $x \in
dom(\sigma_1)$ and $(\sigma(x) \; \overrightarrow{N[\sigma]})
\rhd^* \m \al M_1$ (resp. $\l y M_1$).
\end{lemma}

\begin{proof}
A $\m$-substitution cannot create a $\l$ or a $\m$ (see, for
example,  \cite{DN2}) and thus, the proof is as in lemma
\ref{lambda}.
\end{proof}

\ignore{By induction on $M$. $M$ cannot be of the form
$(\beta \,
M')$ or $\l x \, M'$. If $M$ begins with a $\m$, the
result is
trivial. Otherwise $M=(M_1 \; M_2)$ and, by lemma
\ref{ll2},
$M_1[\sigma] \ras \m\al R$ and $R[\al=_rM_2[\sigma]]
\ras P$. By
the induction hypothesis $M_1\ras \m\al Q$ for some
$Q$ such that
$Q[\sigma] \ras R$ and thus $M \ras \m\al Q[\al=_r
M_2]$. Since
$Q[\al=_r M_2][\sigma] = Q[\sigma][\al=_r M_2[\sigma]]
\ras
R[\al=_r M_2[\sigma]]\ras P$ we are done.
 }

\begin{lemma}\label{mu.untyped}
Assume $M, P, \overrightarrow{Q} \in SN$ and $(M \; P \;
\overrightarrow{Q}) \not\in SN$. Then either ($M \ras \l x M_1$
and $(M_1[x:=P] \;\overrightarrow{Q}) \not \in SN$) or ($M \ras \m
\al M_1$ and $(\m \al M_1[\al =P] \;\overrightarrow{Q}) \not \in
SN$).
\end{lemma}
\begin{proof}
As in lemma \ref{untyped}.
\end{proof}

\begin{lemma}\label{mu.subterm}
Let $M$ be a term and $\sigma$ be a $\l$-substitution.
Assume $M,
\sigma \in SN$ and $M[\sigma] \not\in SN$. Then
$(\sigma(x)
\;\overrightarrow{P[\sigma]}) \not \in SN$ for some
$(x \;
\overrightarrow{P}) \preceq M$ such that
$\overrightarrow{P[\sigma]}\in SN$.
\end{lemma}
\begin{proof}
As in lemma \ref{subterm}.
\end{proof}

\begin{definition}
A $\mu$-substitution $\sigma$ is said to be fair if, for each
$\alpha \in dom(\sigma)$, $\alpha \not \in Fv(\sigma)$ where $x
\in Fv(\sigma)$ (resp. $\beta \in Fv(\sigma)$) means that $x \in
Fv(N)$
  (resp. $\beta \in Fv(N)$) for some $N \in Im(\sigma)$.

\end{definition}

\begin{lemma}\label{goodsub}
Let $\sigma$ be is a fair $\mu$-substitution, $\alpha \in
dom(\sigma)$ and $x \not \in Fv(\sigma)$ (resp. $\beta \not \in
Fv(\sigma)$), then $M[\sigma][x:=\sigma(\alpha)]=
M[x:=\sigma(\alpha)][\sigma]$ (resp. $M[\sigma][\beta
=\sigma(\alpha)]= M[\beta =\sigma(\alpha)][\sigma]$).
\end{lemma}

\begin{proof}
Immediate.
\end{proof}

\begin{lemma}\label{mu.nonsn}
Let $M,N$ be terms and $\sigma$ be a fair $\m$-substitution.
Assume $M[\sigma], N \in SN$ but $(M[\sigma] \; N) \not \in SN$.
Assume moreover that $M[\sigma] \ras \m\al M_1$.  Then, for some
$(\al \; M_2) \preceq M$, we have $(M_2[\sigma'] \; N) \not\in SN$
and $M_2[\sigma'] \in SN$ where $\sigma'=[\sigma + \al =N]$.
\end{lemma}
\begin{proof}
By lemma \ref{mu.lambda},  we know that $M \ras \m\al M'_1$ for
some $M'_1$ such that $M'_1[\sigma] \ras M_1$.  Let $M'$ be a
sub-term of a reduct of $M$ such that $\langle\eta(M'[\sigma]),
cxty(M')\rangle$ is minimum and $M'[\sigma'] \not\in SN$. We show
that $M'=(\al \; M_2)$ and has the desired properties. By
minimality, $M'$ cannot be of the form $\l x P$, $\m \be P$ nor
$(\be \; P)$ for $\be \neq \al$ or $\be \not\in dom(\sigma)$.

If $M'= (P_1 \; P_2)$. By the minimality of $M'$, $P_1[\sigma'],
P_2[\sigma'] \in SN$. Thus, by lemma \ref{mu.lambda} and
\ref{mu.untyped}, $P_1 \ras \l x Q$ (resp. $P_1 \ras \m \be Q$)
such that $Q[\sigma'][x:=P_2[\sigma']]= Q[x:=P_2][\sigma'] \not\in
SN$ (resp. $Q[\sigma'][\be=P_2[\sigma']]=Q[\be=P_2][\sigma']
\not\in SN$) and this contradicts the minimality of $M'$.

If $M'=(\be \; P)$ for some $\be \in dom(\sigma)$. Then
$(P[\sigma'] \; \sigma(\be)) \not\in SN$ and, by the minimality of
$M'$, $P[\sigma'] \in SN$. Thus, by lemmas \ref{mu.lambda},
\ref{mu.untyped} and \ref{goodsub}, $P \ras \l x Q$ (resp. $P \;
\ras \; \m \gamma Q$) such that $Q[\sigma'][x:=\sigma(\be)]=
Q[x:=\sigma(\be)][\sigma'] \not\in SN$ (resp.
$Q[\sigma'][\gamma=\sigma(\be)]=Q[\gamma=\sigma(\be)][\sigma']
\not\in SN$) and this contradicts the minimality of $M'$.

Thus $M'=(\al \; M_2)$ and its minimality implies $M_2[\sigma']
\in SN$.
\end{proof}

\subsection{Proof of the strong normalization}\label{4.1}

We  use the same notations as in section \ref{proof}.

\begin{lemma}\label{mu.sansX}
 Let ${\cal Y} \subseteq {\cal X}$ be such that $H[\{X\}]$ holds for
 each $X \in {\cal Y}$. Then $H[{\cal T}({\cal Y})]$ holds.
\end{lemma}

\begin{proof} Assume that $H[\{X\}]$ holds for
 each $X \in {\cal Y}$. The result is a special case of the following claim.

{\em Claim :} Let $M$ be a term, $U,V$ be types such that $U \in
{\cal T}({\cal Y})$ and $\sigma$ be a $\l$-substitution such that,
for each $x$, $\sigma(x)=N_x[\tau_x]$ where $\tau_x$ is a fair
$\mu$-substitution such that $dom(\tau_x) \cap Fv(M[\sigma]) =
\emptyset$. Assume $\G \v M : V$ and for each $x \in dom(\sigma)$,
$x : U \in \G$.  Assume finally that $M$ and the $N_x[\tau_x]$ are
in $SN$. Then, $M[\sigma] \in SN$.

{\em Proof.} By induction on $\langle lg(U),\eta c(M), \eta c
(\sigma)\rangle$ where $\eta(\sigma)= \;$  $\sum \eta(N_x)$ and
$cxty(\sigma)=  \sum cxty(N_x)$ and, in the sums, each occurrence
of a variable  counts for one. For example, if there are two
occurrences of $x_1$ and three occurrences of $x_2$,
$cxty(\sigma)=2\;cxty(N_1)+3\;cxty(N_2)$. Note that we really mean
$cxty(N_x)$ and not $cxty(N_x[\tau_x])$ and similarly for $\eta$.

The only non trivial case is when  $M = (x\;
Q\;\overrightarrow{O})$ for $x \in dom(\sigma)$. By the {\em IH},
$Q [\sigma] , \overrightarrow{O[\sigma]}\in SN$. It is enough to
show that $(N_x[\tau_x] \; Q[\sigma]) \in SN$ since $M[\sigma]$
can be written as $M'[\sigma']$ where $M'=(z \;
\overrightarrow{O[\sigma]})$ and $\sigma'(z)=(N_x[\tau_x] \;
Q[\sigma])$ and (since the size of the type of $z$ is less than
the one of $U$) the \ihb gives the result. By
 lemma \ref{mu.untyped}, we have two cases to
consider.

\begin{itemize}
 \item $N_x[\tau_x] \ras \l y N_1$. By lemma \ref{mu.lambda}, $N_x
\ras \l y N_2$ and the proof is exactly the same as in lemma
\ref{sansX}.
\item $N_x[\tau_x] \ras \m \al N_1$. By lemma \ref{mu.nonsn}, let
$(\al \; N_2) \preceq N_x$ be such that $N_2[\tau'] \in SN$ and
$R=(N_2[\tau'] \; Q[\sigma]) \not \in SN$ where $\tau'=[\tau_x +
\al =Q[\sigma]]$. But $R$ can be written as $(y \; Q)[\sigma']$
where $\sigma'$ is the same as $\sigma$ except that $\sigma'(y)=
N_2[\tau']$. Note that $(y \; Q)$ is the same as (or less than)
$M$  but one occurrence of $x$ has been replaced by the fresh
variable $y$. The substitution $\tau'$ is fair and $dom(\tau')
\cap Fv((y \; Q)) = \emptyset$. The \ihb gives a contradiction
since $\eta c(\sigma') < \eta c(\sigma)$. Note that the type
condition on $\sigma'$ is satisfied since $N_x$ has type $U$, thus
$\al$ has type $\neg U$ and thus $N_2$ also has type $U$.\qed
\end{itemize}
\end{proof}

\medskip

\noindent {\em For now on, we fix some $i$ and we assume
$H[\{X_j\}]$ for each $j < i$.  Thus, by lemma \ref{mu.sansX}, we
know that $H[{\cal T}'_i]$ holds. It remains to prove $H[\{X_i\}]$
i.e. proposition \ref{mu.avecX}.}

\medskip

\begin{definition}
Let $M$ be a term,  $\sigma=\sigma_1 \cup \sigma_2$ where
$\sigma_1$ (resp. $\sigma_2$) is a $\l$ (resp. $\m$) substitution,
$\G$ be a context and $U$ be a type. Say that
 $(\sigma,\G , M , U)$ is adequate if the following
holds:
\begin{itemize}
\item  $\G \v
M[\sigma] : U$ and $M, \sigma \in SN$.
\item For each $x \in dom(\sigma_1)$, $\G \v \sigma(x) : V_x$ and  $V_x \in {\cal
T}_i^+ $.
\end{itemize}
\end{definition}

Note that nothing is  asked on the types of the
$\m$-variables.

\begin{lemma}\label{mu.prepa}
Let $n,m$ be integers, $\overrightarrow{S}$ be a sequence of terms
and $(\delta, \Delta , P , B)$ be adequate.  Assume that
\begin{enumerate}
\item  $B \in {\cal T}_i^- - {\cal T}'_i$ and $\Delta \v (P
[\delta] \, \overrightarrow{S}) : W$ for some $W$.
\item $\overrightarrow{S} \in SN$, $P \in SN$ and $\eta c(P) < \langle n, m \rangle$.
\item $M[\sigma] \in SN$ for every adequate $(\sigma, \G, M, U)$ such
that $\eta c (M) <  \langle n, m \rangle$.
\end{enumerate}
Then $(P [\delta] \, \overrightarrow{S})  \in SN$.
\end{lemma}

\begin{proof}
By induction on the length of $\overrightarrow{S}$. The proof is
as in lemma \ref{prepa}. The new case is $P[\delta] \rhd^* \m \al
R$ (when
 $\overrightarrow{S} = S_1
\overrightarrow{S_2}$). By lemma \ref{mu.lambda}, we have two
cases to consider.

\begin{itemize}
  \item $P \rhd^* \m \al R'$. We have to show that $Q
= (\m \al R' [\delta + \al = S_1]\ \overrightarrow{S_2}) \in SN$.
By lemma \ref{typed1'},  the properties of $\approx$ and since $B
\in {\cal T}_i^-$, there are types $B_1, B_2$ such that $B \approx
B_1 \rightarrow B_2$ and $\Delta \v \m \al R'[\delta + \al = S_1]
: B_2 $  and $B_2  \in {\cal T}_i^-$. Since $\eta c(R') < \langle
n , m \rangle$ and $([\delta +\al=S_1], \Delta \cup \{ \al : \neg
B_2\}, \mu \al R', B_2)$ is adequate, it follows from (3) that $R'
[\delta + \al= S_1] \in SN$.

- Assume first $B_2 \in {\cal T}'_i$. Since $(z' \,
\overrightarrow{S_2}) \in SN$ and $Q = (z' \, \overrightarrow{S_2})
           [z':= \m\al R' [\delta + \al= S_1]]$, the result
           follows from $H[{\cal T}'_i]$.

- Otherwise, the result follows from the {\em IH} since $([\delta
  +\al=S_1], \Delta \cup \{ \al : \neg B_2\}, \mu \al R', B_2)$ is
adequate and the length of $ \overrightarrow{S_2}$ is less than the
one of $\overrightarrow{S}$.

\item $P \rhd^* (y \, \overrightarrow{T})$ for some $\l$-variable $y
\in dom(\delta)$. As in lemma \ref{prepa}. \qed

\end{itemize}
\end{proof}

\begin{lemma}\label{mu.2}
Assume $(\sigma,  \G , M, A )$ is adequate.  Then $M
[\sigma] \in SN$.
\end{lemma}

\begin{proof}
As in the proof of the lemma \ref{mu.sansX}, we  prove a more
general result.  Assume that, for each $x \in dom(\sigma_1)$,
$\sigma_1(x)=N_x[\tau_x]$ where $\tau_x$ is a fair
$\mu$-substitution such that $dom(\tau_x) \cap Fv(M[\sigma]) =
\emptyset$. We  prove that $M[\sigma] \in SN$.

By induction on $ \eta c(M)$ and, by secondary induction, on $\eta
c(\sigma_1)$ where $ \eta(\sigma_1)$ and $cxty(\sigma_1)$ are
defined as in lemma \ref{mu.sansX}. The proof is as in lemma
\ref{mu.sansX}. The interesting case is $M = (x\;
Q\;\overrightarrow{O})$ for some $x \in dom(\sigma_1)$. The case
when $ N_x[\tau_x] \ras \l y N'$ is as in lemma \ref{X>0}. The new
case is when $N_x[\tau_x] \rhd^* \m \al N'$. This is done as in
lemma \ref{mu.sansX}. Note that, for this point, the type was not
used.
\end{proof}

\begin{proposition}\label{mu.avecX}
Assume $\G , x : X_i \v M : U$ and $\G \v N : X_i$ and $M,N \in
 SN$. Then $M[x := N] \in SN$.
\end{proposition}

\begin{proof}
This follows from lemma \ref{mu.2} since $([x:=N], \G, M, U )$ is
adequate.
\end{proof}

\ignore{
\begin{lemma}\label{mu.crux}
Let $M,N$ be terms in $SN$. Assume $\G, x:U \v M : V$
and $\G \v N
: U$.  Then $M[x := N] \in SN$.
\end{lemma}

\begin{proof}
We prove the same more general result as in lemma \ref{mu.1} using
induction on the same 5-tuple. The proof is exactly the same.  As
there,  the only thing we have to check is that $U$ is an arrow
type. We conclude, as in lemma \ref{crux}, by using  corollary
\ref{mu.avecX}.
\end{proof}

\begin{theorem}
Every typed $\l\m$-term $M$ is in $SN$
\end{theorem}

\begin{proof}
As in theorem \ref{SNF}.
\end{proof}
}

\section{Some applications}\label{applic}

\subsection{Representing more functions}

By using  recursive types, some terms that cannot be typed in the
simply typed $\l$-calculus become typable.  For example, by using
the equation $X \approx (X \f T) \f T$, it is possible to type
terms containing both $(x \; y)$ and $(y \; x)$ as sub-terms. Just
take $x:X$ and $y: X \f T$. By using the equation $X \approx T \f
X$, it is
 possible to apply  an unbounded number of arguments to a term.

 It is thus natural to try to extend  Schwichtenberg's result and to
determine the class of
 functions that are represented in such systems and, in particular, to
see whether or not they
 allow to represent more functions. Note that Doyen \cite{Doy} and
 Fortune \& all \cite{Lei} have given extensions of Schwichtenberg's
 result.

Here is an example of function that cannot be typed (of the good
type) in the simply typed $\l$-calculus.

Let ${\it Nat} = (X \f X) \f (X \f X)$ and ${\it Bool} = Y \f (Y
\f Y)$ where $X,Y$ are  type variables. Let $\tilde{n} = \l f \l x
\; (f \; (f \; ... \; x)\; ...)$ be the church numeral
representing $n$ and ${\bf 0} = \l x \l y \; y$,  ${\bf 1} = \l x
\l y \; x$ be the terms representing {\it false} and {\it true}.
Note that $\tilde{n}$ has type {\it Nat} and ${\bf 0}$, ${\bf 1}$
have type ${\it Bool}$.

The term ${\it Inf} = \l x\l y \;(x \, M \, \l z {\bf 1} \, (y \,
M \, \l z {\bf 0}))$ where $M = \l x \; \l y \; (y \, x)$  has
been introduced by B.Maurey. It is easy to see that, for every
$n,m \in \mathbb{N}$, the term $({\rm Inf} \, \widetilde{m} \,
\widetilde{n})$ reduces to {\bf 1} if $m\leq n$ and to {\bf 0}
otherwise. Krivine has shown in \cite{Kri2} that the type ${\it
Nat} \f  {\it Nat} \f {\it Bool}$ cannot be given to ${\it Inf}$
in system $F$ but, by adding the equation $X \approx (X \f {\it
Bool}) \f {\it Bool}$, it becomes typable. Our example uses the
same ideas.

\bigskip

  Let  $\approx$ be the congruence
 generated by $X \approx  (X \f {\it Bool}) \f
 {\it Bool}$.
For each $n \in \N^*$, let ${\it Inf}_n = \l x \; (x \, M \, \l y
{\bf 1} \, (M^{n-1} \, \l y {\bf 0}))$ where $(M^k \, P)= (M \; (M
\; ... \; (M \; P)))$.

\begin{proposition}\label{p1}
For each $n \in \N^*$ we have $\v {\it Inf}_n : {\it Nat} \f {\it
Bool}$.
\end{proposition}
\begin{proof}
We have $x : X \f {\it Bool}, y : X \v (y \, x) : {\it Bool}$,
then $\v M : (X \f {\it Bool}) \f (X \f {\it Bool})$, thus $\v
(\widetilde{n} \, M) : (X \f {\it Bool}) \f (X \f {\it Bool})$.
But $\v \l y  {\bf 0} : X \f {\it Bool}$, therefore $\v
(\widetilde{n} \, M \, \l y {\bf 0}) : X \f {\it Bool}$.

We have $x : X , y : X \f {\it Bool} \v (y \, x) : {\it Bool}$,
then $\v M : X \f X$, thus $x : {\it Nat} \v (x \, M) : X \f X$.
But $\v \l y  {\bf 1} : (X \f {\it Bool}) \f {\it Bool}$,
therefore $x : {\it Nat} \v (x \, M \, \l y
 {\bf 1}) : X$.

We deduce that $x : {\it Nat} \v ((\widetilde{n} \, M \, \l y {\bf
0}) \, (x \, M \, \l y {\bf 1})) : {\it Bool}$, then $x : {\it
Nat} \v (x \, M \, \l y  {\bf 1} \, (M^{n-1} \, \l y {\bf 0})) :
{\it Bool}$ and thus $\v {\it Inf}_n : {\it Nat} \f {\it Bool}$.
\end{proof}

\begin{proposition}\label{p2}
For each $n \in \N^*$ and $m \in \N$, $({\it Inf}_n \,
\widetilde{m})$ reduces to  {\bf 1} if $m \leq n$ and to {\bf 0}
otherwise.
\end{proposition}
\begin{proof}

 $({\it Inf}_n \, \widetilde{m}) \; \rhd^*  (M^m \, \l y  {\bf
1} \, (M^{n-1} \, \l y {\bf 0}))
 \; \rhd^*  (M^{n-1} \, \l y  {\bf 0} \, (M^{m-1} \,
\l y {\bf 1})) \; \rhd^*$  \\
$(M^{m-1} \, \l y {\bf 1} \, (M^{n-2} \, \l y {\bf 0})) \; \rhd^*
(M^{n-2} \, \l y {\bf 0} \, (M^{m-2} \, \l y {\bf 1})) \; \rhd^*
... $ \\
$\rhd^* {\bf 1 \;\; } {\mbox if }\;\;  m \leq n \;\; {\mbox and}
\;\; {\bf 0} \;\ {\mbox otherwise}.$
\end{proof}

\noindent {\bf Remarks}

Note that for the (usual) simply typed $\l$-calculus we could have
taken for $X$ and $Y$ the same variable but, for propositions
\ref{p1} and \ref{p2}, we cannot assume that $X=Y$  because then
the condition of positivity would not be satisfied. This example
is thus not completely satisfactory and it actually shows   that
the precise meaning of the question ``which functions can be
represented in such systems'' is not so clear.

\subsection{A translation of the $\l\m$-calculus into the
$\l$-calculus}\label{transfo}

The  strong normalization of a typed $\l\m$-calculus can be
deduced from the one of the corresponding typed $\l$-calculus by
using CPS translations. See, for example, \cite{DeGr1} for such a
translation. There is another, somehow simpler, way of doing such
a translation. Add, for each atomic type $X$, a constant $a_X$ of
type $ \neg \neg X \f X$. Using these constants, it is not
difficult to get, for each type $T$, a  $\l$-term $M_T$ (depending
on $T$) such that $ M_T$ has type $ \neg \neg T \f T$. This gives
a translation of the $\l \mu$-calculus into the $\l$-calculus from
which the strong normalization of the $\l\m$-calculus can be
deduced from the one of the  $\l$-calculus. This translation,
quite different from the CPS translations, has been  used by
Krivine  \cite{Kri3} to code the $\l\m$-calculus with second order
types in the $\l {\cal C}$-calculus.

With recursive equations, we do not have to add the constant $a_X$
since we can use the equation $X \approx \neg \neg X$. We give
here, without proof, the translation. We denote by $S_{\approx}$
the simply typed $\l$-calculus where $\approx$ is the congruence
on $\cal{T}$ (where $\cal{A}=\{\bot\}$) generated by $X \approx
\neg \neg X$ for each $X$ and by $S_{\l\m}$ the usual (i.e.
without recursive types) $\l\m$-calculus.

\medskip

\begin{definition}

\begin{enumerate}
\item We define, for each type
$T$, a closed $\l$-term $M_T$ such that  $\v_{\approx} M_T :
\neg\neg T \f T$ as follows. This is done by induction on $T$.
\begin{itemize}
 \item $M_{\bot} = \l x \;(x \; I)$ where $I=\l x \ x$.
 \item If $X \in {\cal X}$, $M_X= I$.
 \item $M_{U\f V} = \l x\l y \; (M_V \; \l z (x \; \; \l t (z \; (t \;
 y))))$
\end{itemize}

\item We define a translation from $S_{\l\m}$ to $S_{\approx}$ as follows.

\begin{itemize}
\item $x^* = x$.
\item $(\l x \, M)^* = \l x \, M^*$.
\item $(M \, N)^* =  (M^* \, N^*)$.
\item $(\mu \al \, M)^* = (M_U \; \l \al \, M^*)$ if $\al$ has
  the type $\neg U$.
\item $(\al \, M)^* =  (\al \, M^*)$.
\end{itemize}
\end{enumerate}
\end{definition}

For a better understanding, in the translation of $\m\al
  M$ and $(\al \, M)$, we have kept the same name to the  variable $\al$  but
  it should be clear that the translated terms are $\l$-terms with
  only on kind of variables.

\begin{lemma}\label{trtype}
If $\G \v_{\l\m} M : U$  then $\G \v_{\approx} M^* : U$.
\end{lemma}

\begin{lemma}\label{trred}
Let $M,N$ be  typed $\l\m$-terms. If $M \rhd N$, then $M^* \rhd^+
N^*$.
\end{lemma}

\begin{proof}
It is enough to check that $(\m\al M \; N)^* \rhd^+ (\m\al
M[\al=N])^*$.
\end{proof}

\begin{theorem}
The strong normalization of $S_{\approx}$    implies the one of
$S_{\l\m}$.
\end{theorem}

\begin{proof}
By lemmas \ref{trtype} and \ref{trred}.
\end{proof}

\noindent {\bf Remark}

 Note that the previous translation cannot be used to show that the
$\l\m$-calculus with recursive types is strongly normalizing since
having two equations (for example $X \approx \neg \neg X$ and $X
\approx F$) is problematic.

\section{Remarks and open questions}\label{conc}

\begin{enumerate}

  \item The proof of the strong normalization of the
  system D of intersection types
\cite{Coppo} is exactly the same as the one for simple types. Is
it possible to extend our proof  to such systems with equations ?
Note that the sort of constraints that must be given on the
equations is not so clear. For example, what does that mean to be
positive in $A \wedge  B$ ? To be positive  both  in $A$ and $B$ ?
in one of them ? It will be interesting to check precisely
because, for example,  it is known that  the system\footnote{This
example appears in a list of open problems of the working group
Gentzen, Utrecht 1993.} given by system D and the equations $X
\approx (Y \f X) \wedge (X \f X)$ and $ Y\approx X \f Y$ is
strongly normalizing (but the proof again is not formalized in
Peano arithmetic) though the positivity condition is violated.
  \item  We could add other typing rules and constructors to
ensure that, intuitively, $X$ represents the {\em least fixed
point} of the equation $X \approx F$. This kind of thing is done,
for example, in $TTR$. What can be said for such systems?
\item There are many  translations from, for example, the
$\l\m$-calculus into the $\l$-calculus that allows to deduce the
strong normalization of the former by the one of the latter. These
CPS transformations differ from the one given in section
\ref{transfo} by the fact that the translation of a term does not
depend on its type. What is the behavior of such translations
 with recursive equations ?
\end{enumerate}

\bigskip

\noindent {\bf Acknowledgments}

We would like to thank P Urzyczyn who has mentioned to us the
question solved here and  has also indicated some errors appearing
in previous versions of our proofs. Thanks also to the referees
and their valuable remarks.

\end{document}